\documentclass[11pt]{article}
\input epsf
\usepackage{amsfonts}
\usepackage{amscd}
\usepackage{amsmath,amssymb,amsthm}
\usepackage{amstext}
\usepackage{amscd}
\usepackage{graphicx}
\usepackage[all]{xy}
\usepackage{pstricks,pst-node,pst-tree}

\theoremstyle{plain}
\newtheorem{lem}{Lemma}[section]
\newtheorem*{main}{Theorem 1}

\newtheorem{theo}[lem]{Theorem}

\newtheorem{coro}[lem]{Corollary}

\theoremstyle{definition}

\newtheorem{definition}[lem]{Definition}
\newtheorem{rems}[lem]{Remarks}
\newtheorem{rem}[lem]{Remark}

\newtheorem{exs}[lem]{Examples}
\newtheorem{conjecture}[lem]{Conjecture}

\newcommand{\R}{\mathbb{R}}

\newcommand{\noi} {\noindent}

\renewcommand{\descriptionlabel}[1]%
       {\hspace{\labelsep}\textsf{#1}}

\DeclareMathOperator{\Int}{Int}

\begin{document}

\title{Smoothing closed gridded surfaces embedded in  $\R^4$\,\,}
\author{ Juan Pablo D\'iaz\thanks{This work was partially supported by CONACyT (M\'exico), FORDECYT 265667}, Gabriela Hinojosa\thanks{This work was partially supported by 
CONACyT (Mexico), 
CB-2009-129939.}, Rogelio Valdez, Alberto Verjovsky\thanks{This work was 
partially supported by
CONACyT (Mexico), CB-2009-129280 and PAPIIT (Universidad
Nacional Aut\'onoma de M\'exico) \#IN106817.}}
\date{February 17, 2016}

\maketitle

\begin{abstract}
\noi We say that a topological $n$-manifold $N$ is a cubical $n$-manifold if it is 
contained in the $n$-skeleton of the canonical cubulation $\cal C$ of $\mathbb{R}^{n+k}$ ($k\geq1$). 
In this paper, we prove that any closed, oriented cubical $2$-manifold has a 
transverse field of 2-planes in the sense of Whitehead and therefore it is smoothable by a 
small ambient isotopy.\end{abstract}

\section{Introduction}

The {\it canonical cubulation} $\cal C$ of $\mathbb{R}^{m}$ is the decomposition into 
hypercubes which are the images of the unit cube
$I^{m}=\{(x_{1},\ldots,x_{m})\,|\,0\leq x_{i}\leq 1\}$ by translations by vectors with 
integer coefficients  \cite{BHV}.

\begin{definition}
Let $N$ be a topological $n$-manifold embedded in $\mathbb{R}^{n+k}$. We say that $N$ is a \emph{cubical manifold} of codimension $k$ if it is contained in the
$n$-skeleton of the canonical cubulation of  $\mathbb R^{n+k}$. When $n=2$ and $k=2$, $N$ is called a \emph{gridded surface} in $\R^4$.
\end{definition}

\noindent Observe that a cubical manifold can be subdivided into simplices to become a PL-manifold.  M. Boege, G. Hinojosa and A. Verjovsky proved in  \cite{BHV} the following theorem.

\begin{theo} Let $N$ be a closed and smooth $n$-dimensional submanifold of $\mathbb{R}^{n+2}$ such that it has a trivial normal bundle.
Then $N$ can be deformed by an ambient isotopy into a cubical manifold.
\end{theo}

\noindent The goal of this paper is to prove a sort of reciprocal theorem for cubical manifolds of dimension two in $\R^4$ (gridded surfaces).  

\begin{main}\label{main}
Any closed, oriented, gridded surface $N$ in $\mathbb{R}^{4}$ is smoothable. More precisely $N$
 admits a transverse field of 2-planes and therefore by a theorem of J. H. C. Whitehead there is an arbitrarily small topological isotopy that moves $N$ 
 onto a smooth surface in $\R^4$ (see  \cite{pugh}, \cite{whitehead}).

 \end{main}

\noindent The relevance of this result is that there exist PL-manifolds which are not smoothable. For instance,
Kervaire (\cite{kervaire}) constructed an example of a PL triangulable closed manifold $M$ of dimension 10
that does not admit any differentiable structure. Therefore the Kervaire manifold cannot be embedded as a codimension two cubical submanifold of $\R^{12}$. Theorem 1 
contrasts with the result of M. H. Freedman stating that every homology 3-sphere embeds topologically in $\R^4$ and if its Rokhlin invariant is one the embedding can never be 
smooth (or even PL) since this would violate Rokhlin signature theorem for spin 4-manifolds \cite{rokhlin}, \cite{freedman}. 

\begin{rem}
Theorem 1 is false for PL-surfaces contained in the 2-skeleton of a PL triangulation of ${\mathbb S}^4$. For instance if $K\subset \mathbb S^3$ is the trefoil knot and we consider the pair $(\mathbb S^3,K)$ as a pair of PL-manifolds with respect to some PL-triangulation of the 3-sphere
and $(\mathbb S^4,\Sigma(K))$ is the suspension of the pair with the canonical suspended PL-structure, then 
there exists a subdivision of a PL-triangulation of ${\mathbb S}^4$ such that $\Sigma(K)$ is embedded in the two skeleton of the triangulation.
In this way we obtain a 2-dimensional PL knot $\Sigma(K)\hookrightarrow {\mathbb S}^4$ which has two points (the vertices of the suspension) where the knot cannot be made PL-locally flat because if it were  PL-locally flat one could take the cyclic branched covering of order 5 of $\mathbb S^4$ along  $\Sigma(K)$ and one can see that this implies that the Poincar\'e sphere would be the boundary of a PL homology disk and this contradicts Rokhlin's theorem.
{\em This means that the cubic structure plays an important role in the theorem}.

\end{rem}

\noi

\section{Transverse fields for gridded surfaces}

Let $M$ be a topological $m$-manifold embedded in $\mathbb{R}^{m+k}$. We say that an affine $k$-plane $T\subset\mathbb{R}^{m+k}$ is 
\emph{transverse} to $M$ at $p\in M$ if
$p\in T$ and locally $M$ is the graph of a Lipschitz function $h:T^{\bot}(r)\rightarrow N$, where $T^{\bot}(r)$ denotes the $k$-disk of radius $r$ at $p$
in the affine plane $T^{\bot}\subset\mathbb{R}^{m+k}$ perpendicular to $T$ at $p$ (see \cite{pugh}). 
The map $H:x\mapsto x+h(x)$ is a \emph{graph chart} for
$M$ at $p$ that it sends $T^{\bot}(r)$ homeomorphically onto a neighborhood of $p$ in $M$. Since $H$ and $H^{-1}$ are Lipschitz, $H$
is a \emph{Lipeomorphism}.̣\\

\noindent Let $G=G(k,m+k)$ denote the \emph{Grassmann of $k$-planes in} $\mathbb{R}^{m+k}$. A continuous map $\psi:M\rightarrow G$ such that each affine plane
$p+\psi(p)$ is transverse to $M$ at $p$ is a \emph{transverse field} for $M$.\\

\noindent In \cite{whitehead} was proved that if $M$ admits a transverse field then it has an ambient smoothing: it is Lipeomorphic to a nearby 
${\cal{C}}^{\infty}$ submanifold of $\mathbb{R}^{m+k}$. We will use this result to prove Theorem 1.\\

\begin{figure}[h]  
\begin{center}
\includegraphics[height=3cm]{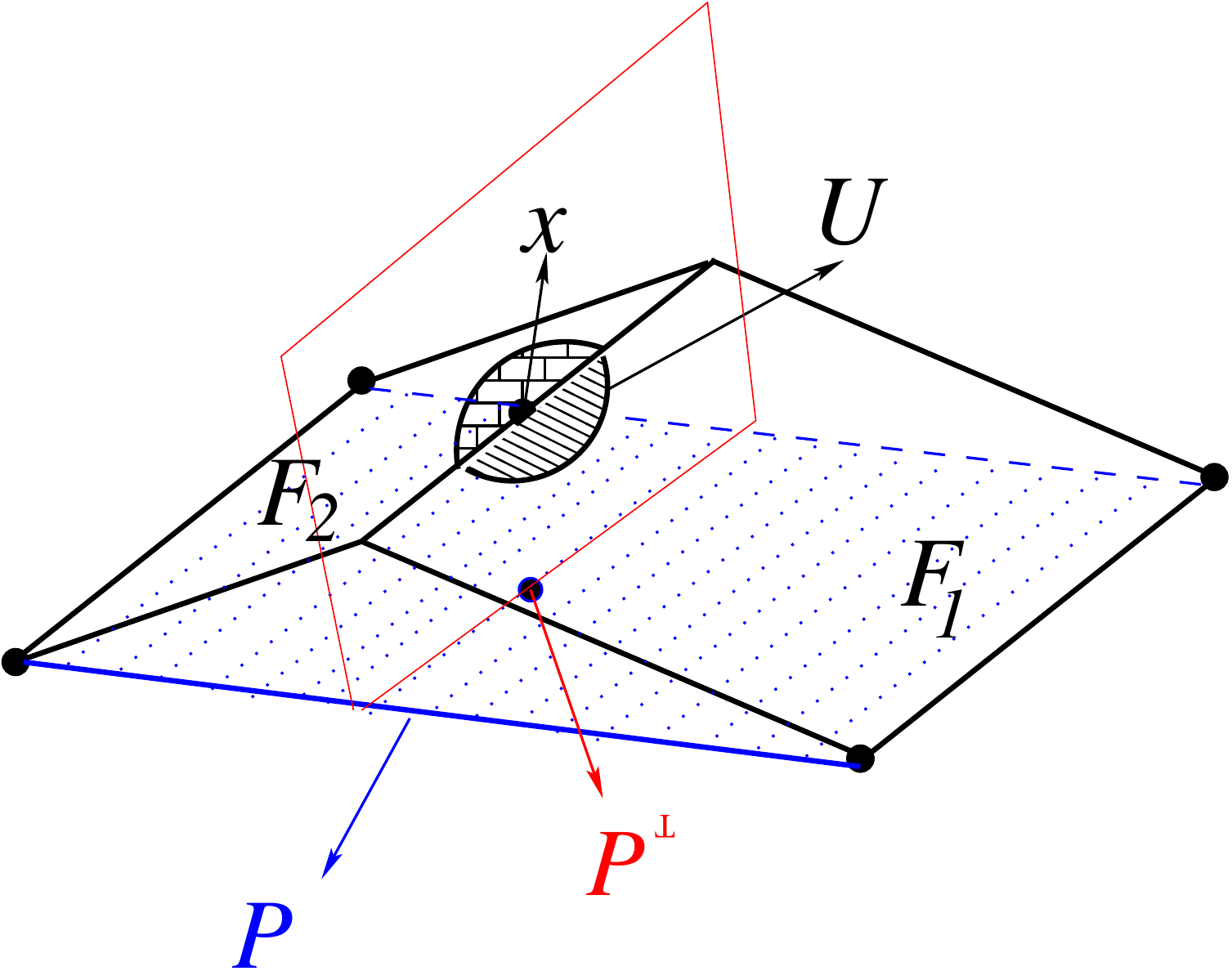}
\end{center}
\caption{\sl The locally transverse vector field ${\text{{\bf v}}}_x$ at $x$.} 
\label{F1}
\end{figure}

\noindent Let $N$ be a closed, oriented, gridded surface  in $\mathbb{R}^{4}$ and consider $x\in N$.  
Notice that if $x$ lies in the interior of some face $F\subset N $,  then we take the neighborhood $U_x$ of $x$ given by $U_x = \Int(F)$. Now consider the plane $P$ parallel to the  support plane of $F$ such that $0\in P$.
Then the map ${\text{{\bf v}}}_x:U_x\rightarrow G(2,4)$ defined as $y\mapsto y+P^{\bot}$, $y\in U_x$, is an example of a locally transverse 2-field to $N$ at $x$, 
where $P^{\bot}$ is the orthogonal plane to $P$.
If the point $x$ lies on the interior of some edge $e$, then 
$x$ belongs exactly to two faces $F_1$ and $F_2$ of $N$, so we can find a neighborhood  $U_x$ of $x$ such that $U_x\subset \Int(F_1\cup F_2)$ and 
$U_x\cap (F_1\cap F_2)=\Int(F_1\cap F_2)$. 
We claim that there exists a locally transverse 2-vector field to $N$ at $x$. In fact, consider
the plane $P$ parallel to the plane generated by the four vertex points of the boundary of $F_1\cup F_2$ (see Figure \ref{F1}) and consider the orthogonal plane $P^{\bot}$ 
to $P$ such that $0\in P$. Then the map ${\text{{\bf v}}}_x:U_x\rightarrow G(2,4)$ given by $y\mapsto y+P^{\bot}$, $y\in U_x$, is an example of a locally transverse 
2-vector field at $x$.

\noindent By the above, the only difficulty lies on the case where $x$ is a vertex of the canonical cubulation $\cal C$. 
 \begin{definition} 
Let $M$ be a topological $n$ manifold. We say that a point $x\in M$ is 
{\it topologically locally flat} or {\it topologically locally tame} if
there exists an open neighborhood $U$ of $x$ such that
there is a homeomorphism of pairs: 
$(U,U\cap M)\sim (\Int(\mathbb{B}^{n+2}),\Int(\mathbb{B}^{n}))$. We say that
$M$ is {\it topologically locally flat}  if all its points are topologically  locally flat. We say that
$M$ is {\it Whitehead locally flat} if it admits a transverse plane field in the sense of Whitehead. If $M$ is
Whitehead  locally flat it is also topologically locally flat. The embedded Poincar\'e sphere in $\mathbb{R}^4$, $\Sigma\hookrightarrow \mathbb{R}^4$ is 
 topologically  locally flat but it is not PL (or Whitehead) locally flat.
\end{definition}

\noindent We will discuss the case where $x$ is a vertex. We can assume without loss of generality and for the sake of simplicity that the vertex $x$ is the origin $0=(0,0,0,0)\in\mathbb{R}^4$ and 
take all the 2-faces of $N$, $F_1,\,F_2,\,\ldots,\,F_j$,  that contain it. Notice that
each edge of $N$ must belong only to two faces  and since
there are only eight edges whose end-point is $0$ it follows that $j\leq 8$. \\

\noindent Consider the unitary canonical vectors on $\mathbb{R}^4$: $e_{\pm 1}=(\pm 1,0,0,0)$, $e_{\pm 2}=(0,\pm 1,0,0)$, 
$e_{\pm 3}=(0,0,\pm 1,0)$ and $e_{\pm 4}=(0,0,0,\pm 1)$. We will use throughout this paper, the following notation for this kind of 2-faces:
$$
F_{u,v}=\{ae_u+be_v\,:\,0\leq a,\,b\leq 1\}
$$ 
where $e_u$ and $e_v$ ($u,v\in\{\pm 1,\,\pm 2,\,\pm 3,\,\pm 4\}$, $|u|\neq |v|$),  denote the corresponding unitary canonical vectors.

\begin{definition}
Let ${\cal {C}}^0$ be the union of all 2-faces $F\in {\cal C}$ such that $0\in F$ and let $N$ be a gridded surface. 
The intersection ${\cal{F}}(N)=N\cap {\cal {C}}^0$ is called the \emph{squared-star} of $N$ at $0$. \\
\end{definition}

\begin{definition}
Let ${\mathbb {F}}$ be the set of all \emph{squared-stars} of all possible gridded surfaces which have a vertex at $0$. \\
\end{definition}

\begin{rem}\label{manifold}
If  the face $F_{a,b}$ belongs to ${\cal{F}}(N)$, since $N$ is a closed 2-manifold, 
there must exist  two 2-faces $F_{a,c},\, F_{b,c}\in {\cal{F}}(N)$. This allows us to describe ${\cal{F}}(N)$ 
as a finite path $\square\,\rightarrow \, \square\, \cdots \rightarrow \, \square$  of consecutive squares \emph{i.e.} squares sharing an edge ($\square$ denotes one of the squares $F_{i,j}\in{\mathcal F}(N))$. 
For instance
$F_{1,2}\,\rightarrow \, F_{1,3}\,\rightarrow \, F_{2,3}$ (see Figure \ref{F2}). Moreover 
\begin{figure}[h]  
\begin{center}
\includegraphics[height=3cm]{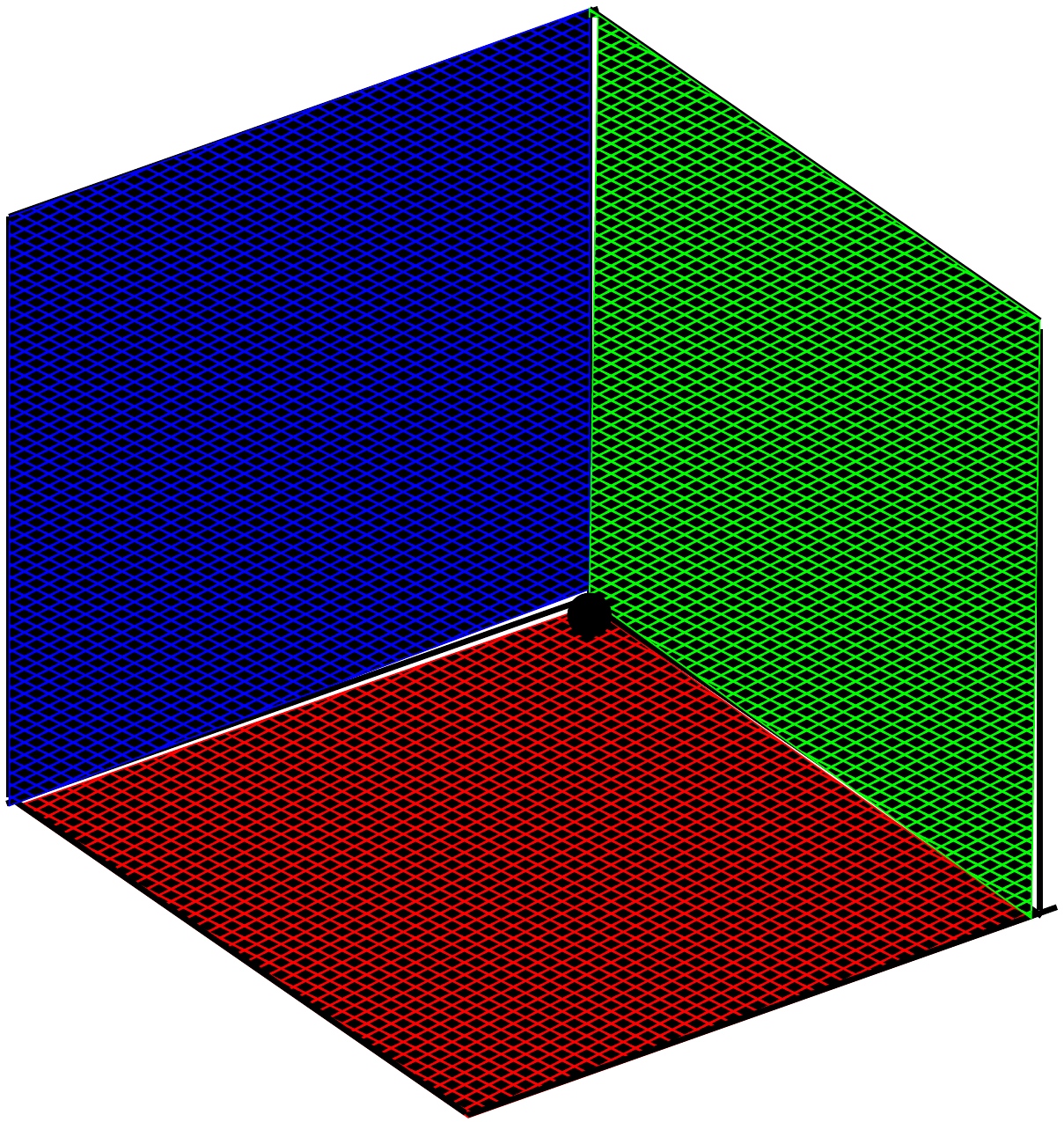}
\end{center}
\caption{\sl squared-star consisting on $3$ squares.} 
\label{F2}
\end{figure}
\end{rem}

\begin{exs}
The number of square faces of a squared-star ${\cal{F}}\in{\mathbb F}$ is any number from 3 to 8 as it is shown in the following examples. 
\begin{enumerate}
\item $F_{1,2}\,\rightarrow \, F_{1,3}\,\rightarrow \, F_{2,3}$ (see Figure \ref{F2}).
\item $F_{1,2}\,\rightarrow \, F_{1,-2}\,\rightarrow \, F_{2,3}\,\rightarrow \, F_{-2,3}$ (see Figure \ref{F3}).
\begin{figure}[h]  
\begin{center}
\includegraphics[height=3cm]{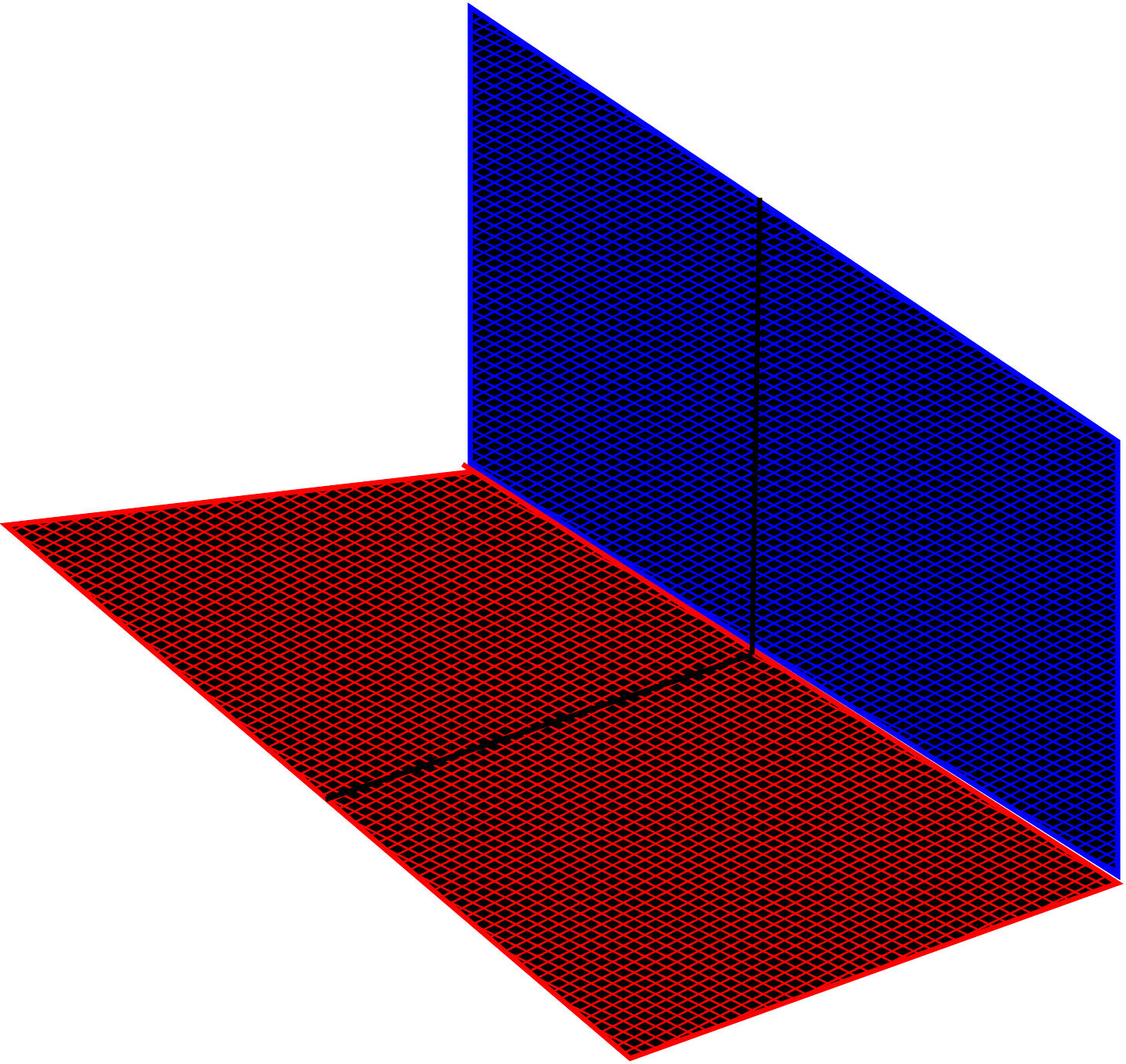}
\end{center}
\caption{\sl squared-star consisting on $4$ squares.} 
\label{F3}
\end{figure}
\item $F_{1,2}\,\rightarrow \, F_{1,-2}\,\rightarrow \, F_{2,3}\,\rightarrow \, F_{-1,3}\,\rightarrow \, F_{-1,-2}$ (see Figure \ref{F4}).
\begin{figure}[h]  
\begin{center}
\includegraphics[height=3cm]{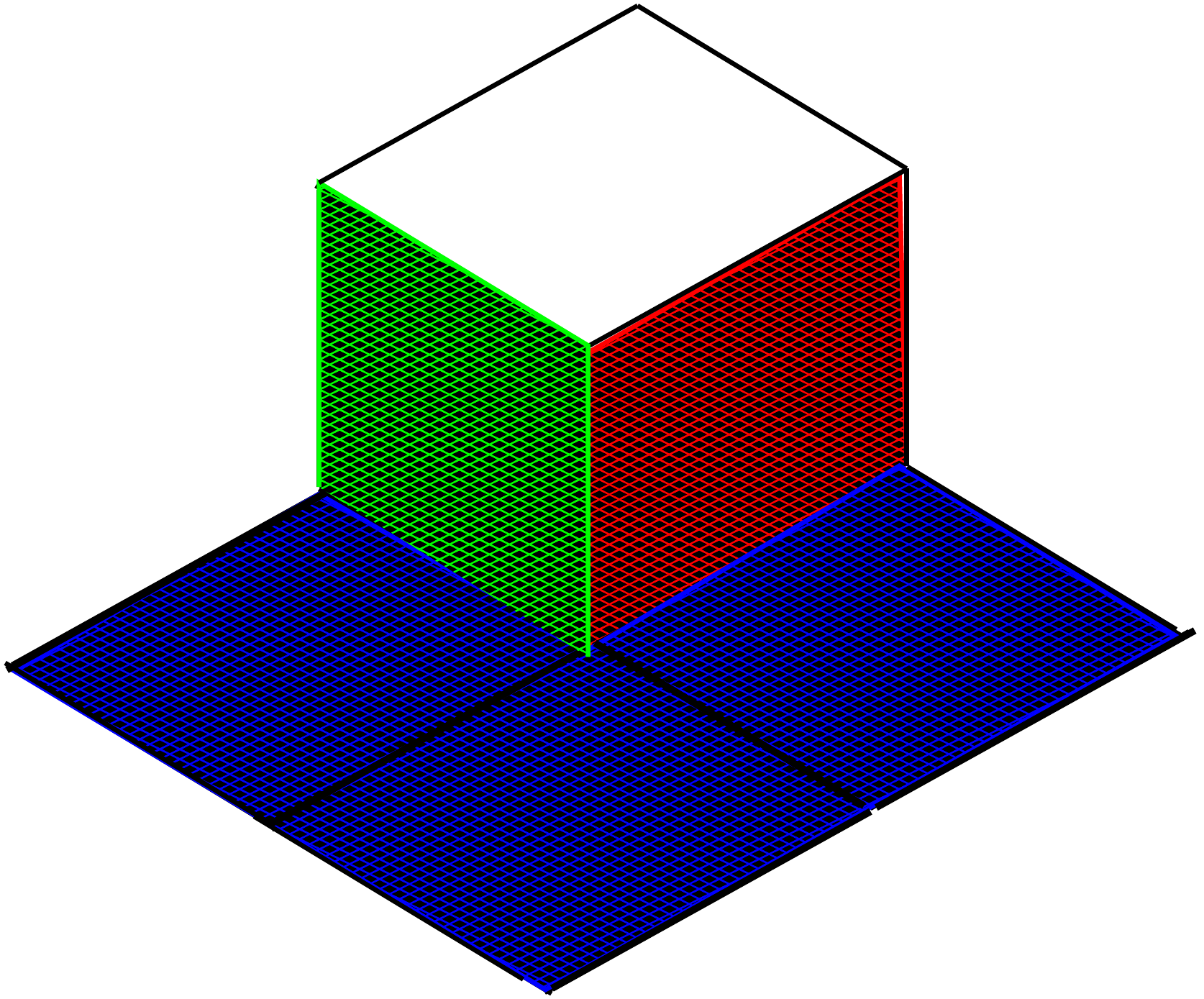}
\end{center}
\caption{\sl squared-star consisting on $5$ squares.} 
\label{F4}
\end{figure}
\item $F_{1,-2}\,\rightarrow \, F_{-1,2}\,\rightarrow \, F_{1,-3}\,\rightarrow \, F_{2,-3}\,\rightarrow \, F_{-2,3}\,\rightarrow \, F_{-1,3}$ (see Figure \ref{F5}).
\begin{figure}[h]  
\begin{center}
\includegraphics[height=3cm]{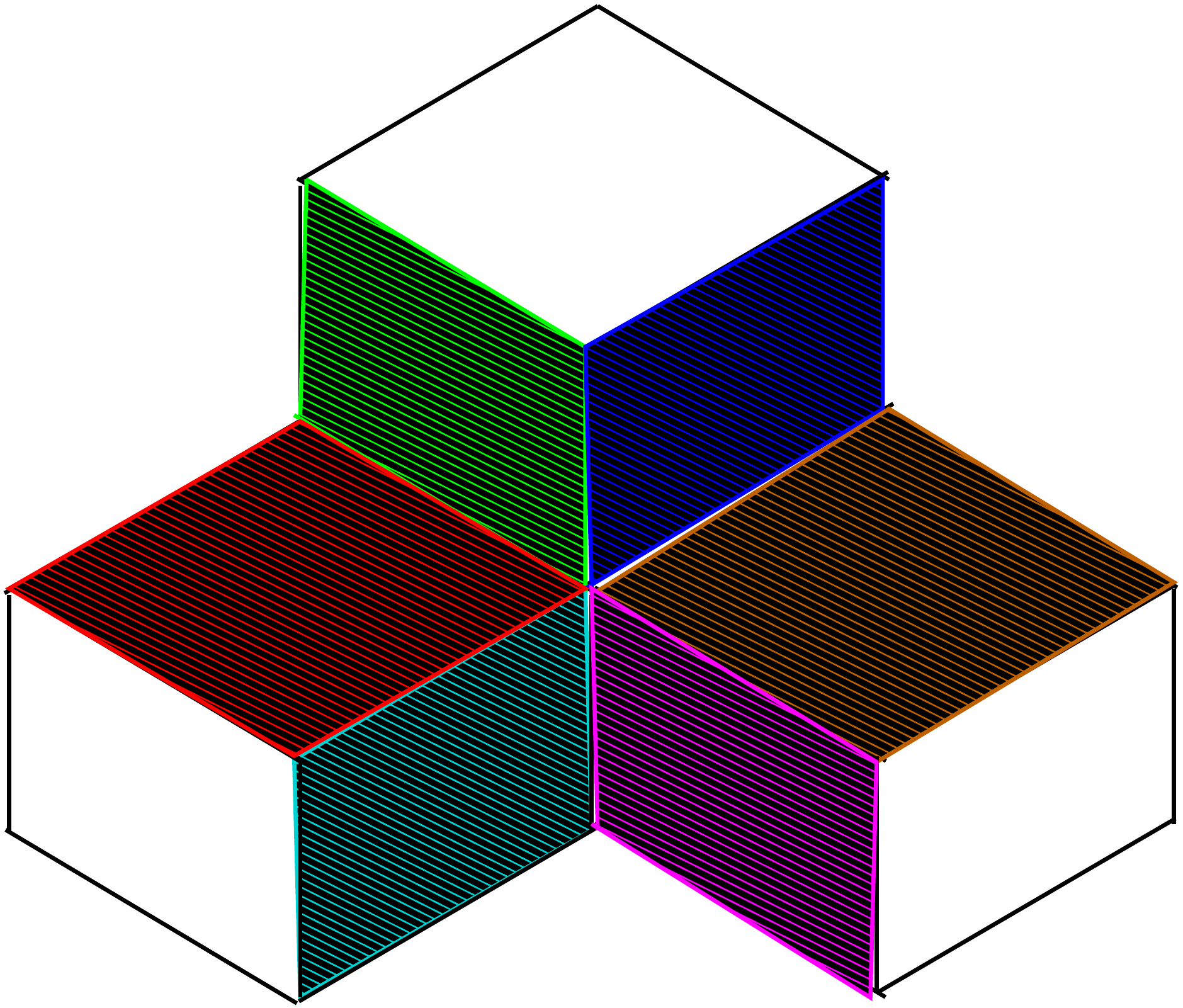}
\end{center}
\caption{\sl squared-star consisting on $6$ squares.} 
\label{F5}
\end{figure}
\item $F_{1,2}\,\rightarrow \, F_{1,-2}\,\rightarrow \, F_{2,3}\,\rightarrow \, 
F_{-1,3}\,\rightarrow \, F_{-2,-3}\,\rightarrow \, F_{-1,4}\,\rightarrow \, 
F_{-3,4}$.
\item $F_{1,2}\,\rightarrow \, F_{1,-2}\,\rightarrow \, F_{-1,2}\,\rightarrow \, 
F_{-1,-3}\,\rightarrow \, F_{-2,-4}\,\rightarrow \, F_{-3,4}\,\rightarrow \, 
F_{3,4}\,\rightarrow \, F_{3,-4}$.
\end{enumerate}
\end{exs}

\noindent Now, we will determine the different squared-stars, up to rotations 
and reflections, in the set $\mathbb{F}$.

\begin{theo} There are 20 different elements in $\mathbb{F}$, up 
to rotations and reflections.
\end{theo}

\noindent The analysis will be made based on the 
number of faces of a squared-star that appear in the same plane. To be able 
to do so, we will think a squared-star as a cycle path of degree 2, as follows.
Consider the complete graph $K_8$ formed with 
the eight vertices labelled $1$, $-1$, $2$, $-2$, $3$, $-3$, $4$, $-4$ 
and edges $F(i,j), \,j\in\{\pm 1, \pm2, \pm 3, \pm 4\}$.
We would like to count the number of closed paths (cycles) of lenght $n$, with 
$3 \leq n \leq 8$, which use the edge $F_{1,2}$, but do not use
the edges $F_{-1,1}$, $F_{-2,2}$, $F_{-3,3}$ or $F_{-4,4}$ and such that every vertex in 
the path has degree $2$.\\

\noindent {\bf Observations.}

\begin{itemize}
\item There are $6$ different planes: generated by the 6 pairs
$e_u$ and $e_v$ ($u,v\in\{1,\, 2,\, 3,\, 4\}$). 

\item Any squared-star ${\cal{F}}$ could have faces in those 6 different 
planes, so the the maximum number 
of planes where there are faces of $F$ is $6$. 

\item If any path has faces $F_{a,b}$ and $F_{c,d}$, with at 
least three of the numbers $a$, $b$, $c$, $d$ having different 
absolute value, then there are 
at least two faces of the squared-star in different planes.

\item There is only one case, where there are 4 faces of a 
squared-star in
the same plane: $F_{1,2} \to F_{2,-1} \to F_{2,-1} \to F_{-1,-2}$.

\item  Except for the previous observation, a squared-star can 
have maximum 3 faces in the same plane.

\item  A squared-star has three faces in the same plane if in the 
path we 
can find three edges of the form: $F_{a,b}$, $F_{b,-a}$, $F_{-a,-b}$.

\item There are two ways that a squared-star can have two faces in 
the same plane: $F_{a,b}\to F_{b,-a}$ or $F_{a,b}\to F_{-a,-b}$.  
\end{itemize}

\noindent We define the signature of a squared-star ${\cal{F}}$ as $(a_1, a_2, \dots, a_k)$, 
the numbers of faces of ${\cal{F}}$ in different planes, where $1 \leq k \leq 6$ and
$a_1 \geq a_2 \geq \dots \geq a_k$. 
For example, the path $F_{1,2} \to F_{2,-3} \to F_{-3,-2} \to
F_{-2,1}$ has signature $(2,2)$, since the path has two faces in the plane 
generated  by $e_1$ and $e_2$ 
and two faces in the plane generated by $e_2$ and $e_3$. 
Except for the one case mentioned in the 
observations, all the numbers in the signature are less than or equal to 3.\\

\noindent The signatures $(1,1,1)$ and $(1,1,1,1)$ are realizable by 
$F_{a,b} \to F_{b,x} \to F_{x,a}$ and $F_{a,b} \to F_{b,c} \to F_{c,d} \to
F_{d,a}$, respectively, if all the symbols have distinct absolute value.

\begin{lem}\label{unos}
Any squared-star cannot have 5 single faces in 5 different planes.
\end{lem}

\noindent {\it Proof}. First let assume that there is
a squared-star with exactly 5 faces, all of them in 5 different planes, that 
is, with signature $(1,1,1,1,1)$. 
The corresponding path will have the form 
$F_{a_1,a_2} \to F_{a_2,a_3} \to F_{a_3,a_4} \to F_{a_4,a_5} \to F_{a_5,a_1}$. Since 
there are $5$ symbols, at least two of them have the same absolute value, 
let say $a_i=-a_j$ with $3 \geq j-i \geq 2$. If $j-i=2$, the path has edges
$F_{a_i,a_{i+1}}$ and $F_{a_{i+1},-a_i}$, which are in the same plane, a 
contradiction. If $j-i=3$, the path will have edges $F_{a_{i-1},a_i}$ and
$F_{-a_i,a_{i+1}}$, a contradiction again.\\

\noindent The three cases left are when a squared-star has, besides the 5 faces 
in different planes, another single face in other plane, two faces in other plane or 
three faces in other plane. But the proof is the same as before, since we will have again in the 
path edges with at least
5 different symbols, and we will arrive to the same contradiction. $\square$\\

\noindent Now, let characterize the squared-stars which have three faces in the 
same plane. Remember that a squared-star has three faces in the same plane, if 
in the path there are three faces of the form $F_{a,b}$, $F_{b,-a}$, $F_{-a,-b}$.

\begin{lem}\label{signature3}
The possible signatures of squared-stars that contains a number $3$ are
$(3,1,1)$, $(3,1,1,1)$, $(3,2,1,1)$, $(3,3,1,1)$.
\end{lem}

\noindent {\it Proof.} Let ${\cal{F}}$ be a squared-star that has 3 appearing in its 
signature, that is, contains three faces of
the form $F_{a,b}$, $F_{b,-a}$, $F_{-a,-b}$. 
We will proceed depending on 
the number of faces in ${\cal{F}}$. It is clear that ${\cal{F}}$ has $n$ faces with
$n \geq 4$. \\

\noindent If ${\cal{F}}$ has 4 faces, they are $F_{a,b}$, $F_{b,-a}$, $F_{-a,-b}$ and
$F_{-b,x}$, then it follows that $x=a$ and the signature of ${\cal{F}}$ is $(4)$, 
a contradiction.\\

\noindent Assume that ${\cal{F}}$ has $5$ faces, that is, 
$F_{a,b}\to F_{b,-a}\to F_{-a,-b}\to F_{-b,t}\to F_{t,a}$. Since
$|t| \neq |a|, |b|$, then the last two faces are in different planes, which 
are not the plane determined by $a$ and $b$, then the signature of ${\cal{F}}$ is
$(3,1,1)$.\\

\noindent If ${\cal{F}}$ has $6$ faces, they have the form
$F_{a,b}\to F_{b,-a}\to F_{-a,-b}\to F_{-b,t}\to F_{t,s}\to F_{s,a}$. Once again,
$|t|, |s| \neq |a|, |b|$ and $|t| \neq |s|$, then it is clear that 
the part $F_{-b,t}\to F_{t,s}\to F_{s,a}$ of the path has signature
$(1,1,1)$. Hence ${\cal{F}}$ has signature $(3,1,1,1)$.\\

\noindent If ${\cal{F}}$ has 7 faces, then the path is
$F_{a,b}\to F_{b,-a}\to F_{-a,-b}\to F_{-b,t}\to F_{t,s}\to F_{s,u}\to F_{u,a}$. 
Then neither $t,s,u$ are equal to $\pm a$ or $\pm b$ and $|t| \neq |s|$. 
Removing $\pm a$ and 
$\pm b$ from the set $\{\pm 1, \pm 2, \pm 3, \pm 4\}$, we obtain that
$t,s, u$ belong to a set of the form $\{\pm p, \pm q\}$. It follows that 
$t=-u$. Then the path becomes
$F_{a,b}\to F_{b,-a}\to F_{-a,-b}\to F_{-b,-u}\to F_{-u,s}\to F_{s,u}\to F_{u,a}$,
which has signature $(3, 2, 1, 1)$.\\
 
\noindent Finally, if ${\cal{F}}$ has 8 faces, it is
$F_{a,b}\to F_{b,-a}\to F_{-a,-b}\to F_{-b,t}\to F_{t,s}\to F_{s,u}\to F_{u,w}\to F_{w,a}$.
 As in the previous case, $t,s,u, w \in \{\pm p, \pm q\}$, then the path is  
$F_{a,b}\to F_{b,-a}\to F_{-a,-b}\to F_{-b,t}\to F_{t,s}\to F_{s,-t}\to F_{-t,-s}\to F_{-s,a}$, with signature $(3,3,1,1)$. $\square$\\

\noindent As we mentioned before, there are two ways that a 
squared-star can have two faces in 
the same plane: $F_{a,b}$, $F_{b,-a}$ or $F_{a,b}$, $F_{-a,-b}$, to the first one 
we will asign the number 2 in the signature and to the second one the number $\bar{2}$.
Now let us find a similar result as Lemma \ref{signature3}, but when $2$ appears
in the signature but 3 does not. \\

\noindent Again let ${\cal{F}}$ be a squared-star with $n$ faces, in whose signature appears 2 but 
not 3. If ${\cal{F}}$ has faces $F_{-a,b}$, $F_{b,-a}$, then it must have
another face $F_{-a,t}$, with $t \neq a$, then $n \geq 4$.
Also let $N=\{\pm 1, \pm 2, \pm 3, \pm 4\} \backslash 
\{\pm a, \pm b\}=\{\pm p, \pm q\}$.\\

\noindent If $n=4$, then ${\cal{F}}$ is $F_{a,b}\to F_{b,-a}\to F_{-a,t}\to F_{t,a}$,
which has signature $(2,2)$.\\

\noindent If $n=5$, then ${\cal{F}}$ is  $F_{a,b}\to F_{b,-a}\to F_{-a,t}\to F_{t,s}\to F_{s,a}$, with $t, s \neq -b$, otherwise the signature of ${\cal{F}}$ has a number 3.
Then $|s|, |t| \neq |a|, |b|$, which implies that the signature is 
$(2,1,1,1)$.\\

\noindent The case $n=6$ gives the path
 $F_{a,b}\to F_{b,-a}\to F_{-a,t}\to F_{t,s}\to F_{s,u}\to F_{u,a}$, where
$t,u \in N$. There are three cases. First
if $s=-b$, the path is  
$F_{a,b}\to F_{b,-a}\to F_{-a,t}\to F_{t,-b}\to F_{-b,u}\to F_{u,a}$, and then
depending whether $|t|=|u|$ or not, we get two possible signatures
$(2,2,\bar{2})$ or $(2,1,1,1,1)$.\\

\noindent If $s \neq -b$, then $t,s,u \in N$, and it follows that the path 
has the form
$F_{a,b}\to F_{b,-a}\to F_{-a,\pm p}\to F_{\pm p, \pm q}\to F_{\pm q, \mp p}$, 
$F_{\mp p,a}$, which have signature $(2,2,\bar{2})$.\\
 
\noindent If ${\cal{F}}$ has 7 faces, the path is 
 $F_{a,b}\to F_{b,-a}\to F_{-a,t}\to F_{t,s}\to F_{s,u}\to F_{u,w}\to F_{w,a}$,
where $t, w \in N$. If $s, u \neq -b$, then the path is
 $F_{a,b}\to F_{b,-a}\to F_{-a,\pm p}\to F_{\pm p, \pm q}\to F_{\pm q,\mp p}\to F_{\mp p, \mp q}\to F_{\mp q, a}$; that is, with signature
$(3,2,1,1)$, a contradiction.\\

\noindent In a similar way, we obtain for $s=-b$ or $u=-b$, 
signatures $(2,2,1,1,1)$ and $(2,\bar{2},1,1,1)$.  \\

\noindent Now, let ${\cal{F}}$ be a squared-star with $n$ faces, in whose signature appears $\bar{2}$ but 
not 3. If ${\cal{F}}$ has edges $F_{a,b}$, $F_{-b,-a}$, then it must have
at least another 4 edges  $F_{b,x}$, $F_{y,a}$, $F_{s,-b}$ and $F_{-a,t}$.
That is, the paths of this kind with 6 edges have two possible forms
$$F_{a,b}\to F_{b,x}\to F_{s,-b}\to F_{-b,-a}\to F_{-a,t}\to F_{y,a}$$
or 
$$F_{a,b}\to F_{b,x}\to F_{s,-a}\to F_{-a,-b}\to F_{-b,t}\to F_{y,a},$$
where $x=s$ and $t=y$. The cases $t=y=-x$ and $t=y \neq -x$, in both cases,
give the signatures $(2,2,\bar{2})$, $(\bar{2}, \bar{2},\bar{2})$ and
$(\bar{2},1,1,1,1)$.\\

\noindent For the analysis we made for the case $n=7$, we know that there
are not signatures of the form $(2,2,2,1)$. Let us prove that the
signatures for $n=7$, with at least two numbers 2, are $(2,\bar{2},1,1,1)$
and $(2,2,1,1,1)$. For that, now it is enough to analyze the case
with at least one $\bar{2}$. We have seen that these paths have the form 
$F_{a,b}\to F_{b,x}\to F_{s,-b}\to F_{-b,-a}\to F_{-a,t}\to F_{y,a}$
or $F_{a,b}\to F_{b,x}\to F_{s,-a}\to F_{-a,-b}\to F_{-b,t}\to F_{y,a}$. We will work with the 
first path, and similarly it will follows for the second one.\\

\noindent If $t=y$, we can add an edge connecting $F_{b,x}$ and $F_{s,-b}$,
obtaining the path
$$F_{a,b} \to F_{b,x}\to F_{x,s} \to F_{s,-b} \to F_{-b,-a} \to F_{-a,t} \to F_{t,a},$$
where $x,s \in \{\pm 1, \pm 2, \pm 3, \pm 4\} \backslash \{ \pm a, \pm b, t\}$.
It follows that the path has signature $(2,\bar{2},1,1,1)$. In the same way, if $x=s$
and adding an edge between $F_{-a,t}$ and $F_{y,a}$, we obtain the same signature  
$(2, \bar{2},1,1,1\}$. Therefore there are not signatures for $n=7$,
with three numbers 2 (or $\bar{2}$) and one number 1.\\

\noindent Similarly there are not signatures for $n=8$ with three numbers 2 and
two numbers 1. That is, the paths with signatures $(2,2,2,1)$ and 
$(2,2,2,1,1)$, (here $2$ can be $\bar{2}$ too) are not realizable. 
Now, let us state the result, that can be proved with the previous analysis.

\begin{lem} \label{signature2}
The signatures for paths of 
length $n \leq 7$, with at least one number 2 and not numbers 3 are
$$
(2,2), (2,1,1,1), (2,2,\bar{2}), (\bar{2}, \bar{2}, \bar{2}), 
(2,2,1,1,1),
$$ 
$$
(2, \bar{2}, 1,1,1), (2,1,1,1,1), (\bar{2},1,1,1,1). \,\,\square
$$
\end{lem}

\noindent We are ready to  determine the different squared-stars, up to 
rotations and reflections, in the set $\mathbb{F}$. As mentioned before, 
the analysis will be made based in the 
number of faces of a squared-star that appear in the same plane.

\noindent Given a number $3 \leq m \leq 8$, we will find the 
different ways to write $m$ as sum of numbers less than or equal to 3, 
and then find if those numbers in the sums form possible signatures.
For example, two ways to write $5$ are $2+2+1$ and $3+2$, but $(2,2,1)$ and
$(3,2)$ are not realizable signatures by Lemmas \ref{signature3} 
and \ref{signature2}.

\begin{itemize}
\item For $n=3$, the ways to write 3 as sums are
$1+1+1$, $2+1$, $3$. The only possible signature for $n=3$ is $(1,1,1)$, which
is realizable by $F_{a,b}$, $F_{b,x}$, $F_{x,a}$, since $x \neq a, b, -b$.
All these squared-stars are images under rotations or reflectios of each other.
Then there is only one path of length 3, up to rotation and reflection.

\item For $n=4$, we have $4=4=3+1=2+1+1=2+2=1+1+1+1$, for which only 
the signatures $(4)$, $(2,2)$ and $(1,1,1,1)$ are realizable, by 
Lemmas \ref{signature3} and \ref{signature2}.

\item For $n=5$, we obtain two possible signatures $(3,1,1)$ and 
$(2,1,1,1)$.

\item For $n=6$, the possible signatures are 
$(3,1,1,1)$, $(2,1,1,1,1)$, $(\bar{2},1,1,1,1)$, $(2,2,\bar{2})$, 
$(\bar{2}, \bar{2}, \bar{2})$.

\item For $n=7$, the signatures are
$(3,2,1,1)$, $(2,2,1,1,1)$, $(2,\bar{2},1,1,1)$, $(\bar{2},\bar{2},1,1,1)$.

\item For $n=8$, the signatures are $(3,3,1,1)$, $(2,2,1,1,1,1)$, 
$(2,\bar{2},1,1,1,1)$, $(2,2,\bar{2},\bar{2})$, $(\bar{2}, \bar{2},
\bar{2}, \bar{2})$.
\end{itemize}

\noindent For the case $n=8$, by Lemmas \ref{unos}, \ref{signature3} and 
\ref{signature2}, those are the only possible signatures and it is easy to 
see that they are realizable; next we will give examples of all of them. $\square$\\

\begin{exs}
We will exhibit a path representative for each signature given above.
\begin{itemize}
\item {\bf{$n=3$}.} $F_{1,2}\,\rightarrow \, F_{1,3}\,\rightarrow \, F_{2,3}$.
\item {\bf{$n=4$}.} Signature
\begin{description} 
\item[$(4)$:] $F_{1,2}\,\rightarrow \, F_{1,-2}\,\rightarrow \, F_{-1,2}\,\rightarrow \, F_{-1,-2}$ 
\item[$(2,2)$:] $F_{1,2}\,\rightarrow \, F_{1,-2}\,\rightarrow \, F_{2,3}\,\rightarrow \, F_{-2,3}$ 
\item[$(1,1,1,1)$:]  $F_{1,2}\,\rightarrow \, F_{2,3}\,\rightarrow \, F_{3,4}\,\rightarrow \, F_{1,4}$
\end{description}
\item {\bf{$n=5$}.} Signature
\begin{description} 
\item[$(3,1,1)$:] $F_{1,2}\,\rightarrow \, F_{-1,2}\,\rightarrow \, F_{-1,3}\,\rightarrow \, F_{-2,3}\,\rightarrow \,F_{1,-2}$
\item[$(2,1,1,1)$:] $F_{1,2}\,\rightarrow \, F_{-1,2}\,\rightarrow \, F_{1,3}\,\rightarrow \, F_{3,4}\,\rightarrow \, F_{-1,4}$
\end{description}

\item {\bf{$n=6$}.} Signature
\begin{description} 
\item[$(3,1,1,1)$:] $F_{1,2}\,\rightarrow \, F_{-1,2}\,\rightarrow \, F_{-1,3}\,\rightarrow \, F_{-2,4}\,\rightarrow\, F_{3,4}\,\rightarrow\, F_{1,-2}$
\item[$(2,1,1,1,1)$:] $F_{1,2}\,\rightarrow \, F_{-1,2}\,\rightarrow \, F_{-1,3}\,\rightarrow \, F_{-2,3}\,\rightarrow\,F_{-2,4}\,\rightarrow\,F_{1,4}$
\item[$(\bar{2},1,1,1,1)$:] $F_{1,2}\,\rightarrow \, F_{-1,-2}\,\rightarrow \, F_{1,3}\,\rightarrow \, F_{-2,3}\,\rightarrow\,F_{2,4}\,\rightarrow\,F_{1,4}$
\item[$(2,2,\bar{2})$:] $F_{1,2}\,\rightarrow \, F_{-1,2}\,\rightarrow \, F_{1,3}\,\rightarrow \, F_{-1,-3}\,\rightarrow\,F_{-2,3}\,\rightarrow\,F_{-2,3}$
\item[$(\bar{2}, \bar{2}, \bar{2})$:] $F_{1,2}\,\rightarrow \, F_{-1,-2}\,\rightarrow \, F_{1,3}\,\rightarrow \, F_{-2,3}\,\rightarrow\,F_{2,-3}\,\rightarrow\,F_{-1,-3}$
\end{description}

\item {\bf{$n=7$}.} Signature
\begin{description} 
\item[$(3,2,1,1)$:] $F_{1,2}\,\rightarrow \, F_{-1,2}\,\rightarrow \, F_{-1,4}\,\rightarrow \, F_{3,-4}\,\rightarrow\, F_{3,4}\,\rightarrow\, F_{-2,-4}\,\rightarrow\,F_{1,-2}$
\item[$(2,2,1,1,1)$:] $F_{1,2}\,\rightarrow \, F_{-1,2}\,\rightarrow \, F_{-1,4}\,\rightarrow \, F_{3,4}\,\rightarrow\, F_{-2,3}\,\rightarrow\, F_{-2,-3}\,\rightarrow\,F_{1,-3}$
\item[$(2,\bar{2},1,1,1)$:] $F_{1,2}\,\rightarrow \, F_{-1,2}\,\rightarrow \, F_{1,-3}\,\rightarrow \, F_{-1,3}\,\rightarrow\, F_{-2,3}\,\rightarrow\, F_{-2,4}\,\rightarrow\,F_{-3,4}$
\item[$(\bar{2}, \bar{2}, 1,1,1)$:] $F_{1,2}\,\rightarrow \, F_{-1,-2}\,\rightarrow \, F_{1,-3}\,\rightarrow \, F_{-1,3}\,\rightarrow\, F_{-2,-3}\,\rightarrow\, F_{2,4}\rightarrow\,F_{3,4}$
\end{description}

\item {\bf{$n=8$}.} Signature 
\begin{description}
\item[$(3,3,1,1)$:] $F_{1,2}\,\rightarrow \, F_{-1,2}\,\rightarrow \, F_{-1,4}\,\rightarrow \, F_{3,-4}\,\rightarrow\, F_{3,4}\,\rightarrow\, F_{-3,-4}\,\rightarrow\,F_{-2,-3}\,\rightarrow\,F_{1,-2}$
\item[$(2,2,1,1,1,1)$:] $F_{1,2}\,\rightarrow \, F_{2,-4}\,\rightarrow \, F_{-1,4}\,\rightarrow \, F_{-1,-4}\,\rightarrow\, F_{3,4}\,\rightarrow\, F_{-2,3}\,\rightarrow\,F_{-2,-3}\,\rightarrow\,F_{1,-3}$
\item[$(2,\bar{2},1,1,1,1)$:] $F_{1,-4}\,\rightarrow \, F_{2,-4}\,\rightarrow \, F_{-1,2}\,\rightarrow \, F_{-1,4}\,\rightarrow\, F_{3,4}\,\rightarrow\, F_{-2,3}\,\rightarrow\,F_{-2,-3}\,\rightarrow\,F_{1,-3}$
\item[$(2,2,\bar{2},\bar{2})$:] $F_{1,2}\,\rightarrow \, F_{-1,2}\,\rightarrow \, F_{-1,4}\,\rightarrow \, F_{3,4}\,\rightarrow\, F_{-2,3}\,\rightarrow\, F_{-2,-3}\,\rightarrow\,F_{1,-4}\,\rightarrow\,F_{-3,-4}$
\item[$(\bar{2}, \bar{2}, \bar{2},\bar{2})$:] $F_{1,2}\,\rightarrow \, F_{-1,-2}\,\rightarrow \, F_{1,3}\,\rightarrow \, F_{-1,-3}\,\rightarrow\, F_{2,4}\,\rightarrow\, F_{-3,4}\,\rightarrow\,F_{-2,-4}\,\rightarrow\,F_{3,-4}$
\end{description}
\end{itemize}
\end{exs}

\begin{rem}
For a higher dimensional cubical $n$-manifold $N\subset \mathbb{R}^{n+2} $, with $n+2>4$, we can also define the notion of 
cubical-star of a vertex (\emph{i.e.} it is the union of all the $n$-dimensional cubes in $N$ which contain the vertex). The set ${\mathbb {F}}$ of 
all \emph{cubical-stars} of all possible cubical knots which have a vertex $0$ becomes extremely complicated. However we have the following conjecture: 
\end{rem}
\begin{conjecture}
Any closed, oriented, cubical  $n$-manifold $N$ in $\mathbb{R}^{n+2}$, $n>2$, is smoothable. More precisely $N$
 admits a transverse field of 2-planes and therefore by a theorem of J. H. C. Whitehead there is an arbitrarily small topological isotopy that moves $N$ 
 onto a smooth manifold in $\R^{n+2}$ (see  \cite{pugh}, \cite{whitehead}).
\end{conjecture}

\section{Smoothing cubulated closed 2-manifolds}

As we can see from the previous section the combinatorial description at vertex points of gridded surfaces is very
complicated. Next, we will study squared-stars from the topological and differentiable point of view.

\noindent Let us remember that we are considering a {\it gridded surface} $N\subset\mathbb{R}^4$, that is, $N$ is contained in the scaffolding $\cal{S}$ of the 
canonical cubulation $\cal{C}$ of $\mathbb{R}^{4}$ and $0$ belongs to $N$.

\begin{definition}
We define the {\emph{squared-link}} of $0$,  as the boundary of its squared-star ${\cal{F}}(N)$ and it is denoted
by $slk(N)$.
\end{definition}

\begin{rems}\label{slk}
\begin{enumerate}
\item Let $D$ be the union of all the hypercubes $Q_i\in{\cal{C}}$ such that $0\in Q_i$.
Notice that $i=1,\,2,\,\dots,\,16$. Let $N$ be a gridded-surface such that $O\in N$. Then its squared-star 
${\cal{F}}(N)$ is contained in $D$ and the squared-link $slk(N)$ is contained in the boundary of $D$. 
\item Remember that  ${\cal{F}}(N)=\cup_{i=1}^{n} F_i$, where $F_i$ is a square (2-face).  Since 
$N$ is a topological manifold, we have that only two edges $e_{i_1}$ and $e_{i_2}$ of $F_i$ belong to 
$slk(N)$ (see Remark \ref{manifold}). Hence $slk(N)$ consists on the union of an even number of edges. More precisely, this number is smaller or equal to sixteen.
\end{enumerate}
\end{rems}

\begin{theo}\label{unknotted}
Let $N$ be a gridded surface such that $0\in N$. Then  $slk(N)$ is an unknotted 
simple closed curve.
\end{theo}
\noindent{\it Proof.} Let ${\cal{F}}(N)=\cup_{i=1}^{n}F_i$ be the squared-star of $N$. We will prove that ${\cal{F}}(N)$ is isotopic to the square $[-1,1]^2$ on 
some canonical $ab$-plane; thas is  ${\cal{F}}(N)$ is equivalent to $F_{a,b}\,\rightarrow \, F_{a,-b}\,\rightarrow \, F_{-a,b}\,\rightarrow \, F_{-a,-b}$. 
We will do it using induction over the number of faces $n$ of ${\cal{F}}(N)$. \\

\noindent If $\cal{F}(N)$ consists of $3$ faces  
$F_{a,b}\rightarrow F_{a,c}\rightarrow F_{b,c}$ (see Figure \ref{F2}), we 
will apply the following claim.\\

\noindent{\bf Claim.} Let  ${\cal{F}}(N)=\cup_{i=1}^{n} F_i$ be a squared-star such that $F_1=F_{a,b}$ and $F_2=F_{a,c}$ ($|b|\neq |c|$). 
Then ${\cal{F}}(N)$ is isotopic to the squared-star ${\cal{F}}(N_1)=F_{b,c}\cup (\cup_{i=3}^{n} F_i)$.\\

\begin{figure}[h]  
\begin{center}
\includegraphics[height=3cm]{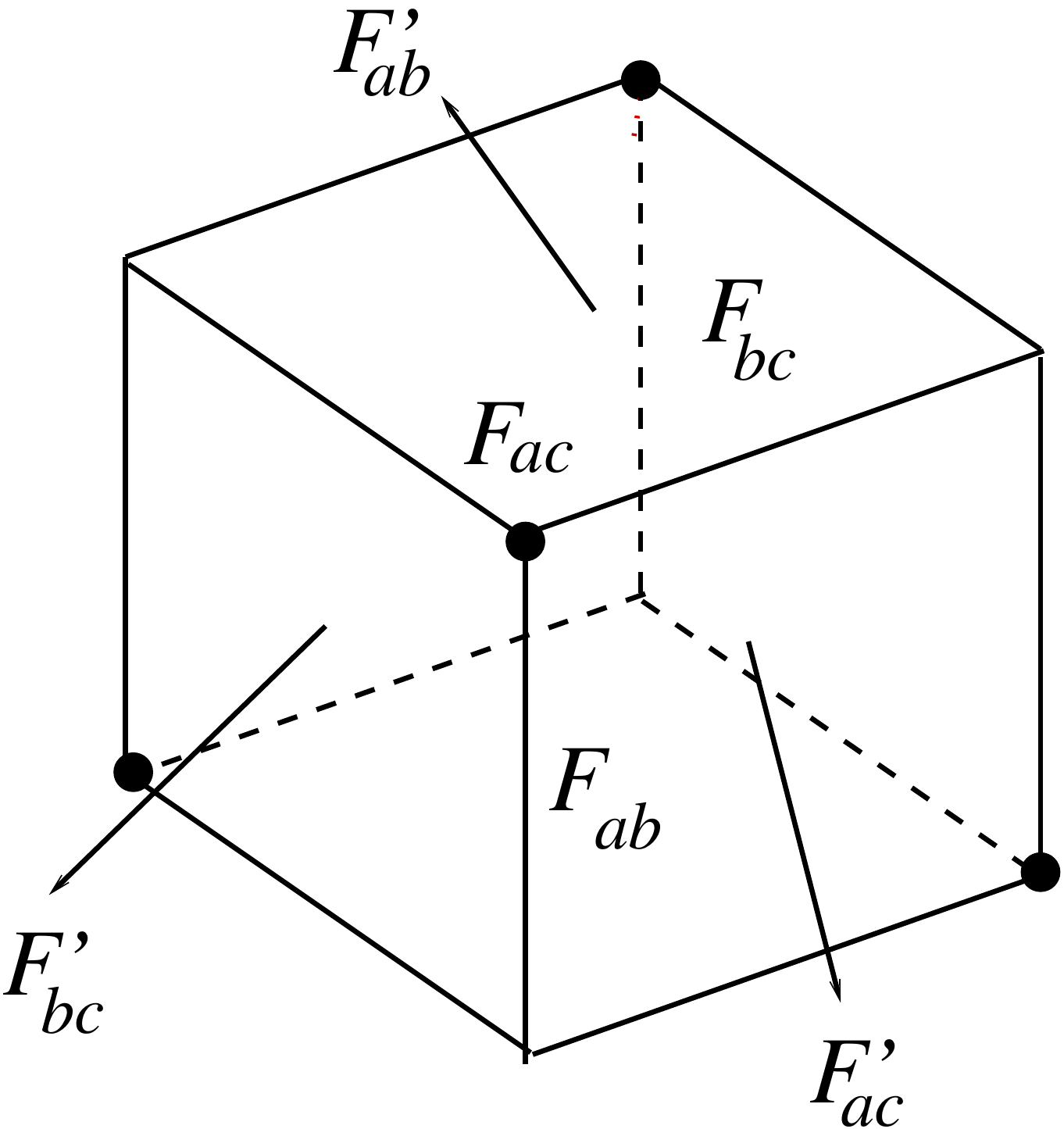}
\end{center}
\caption{\sl All faces of the cube $Q$.} 
\label{F6}
\end{figure}

\noindent{\it Proof of Claim.} Notice that the faces $F_{a,b}$ and $F_{a,c}$  share a common edge $e_a$. Consider the cube $Q$ generated by the canonical vectors 
$e_a$, $e_b$ and $e_c$. Then $Q$ is a 3-face contained in the cubulation ${\cal{C}}$ and $Q$ possesses the 2-faces $F_{a,b}$, $F_{a,c}$, $F_{b,c}$, 
$F'_{a,b}= F_{a,b}+e_c$, $F'_{a,c}= F_{a,c}+e_b$ and $F'_{b,c}= F_{b,c}+e_a$ (see Figure \ref{F6}). Notice that the faces $F'_{a,b}$, $F'_{a,c}$ and $F'_{b,c}$ do not
contain the vertex $0$.

\begin{figure}[h]  
\begin{center}
\includegraphics[height=3cm]{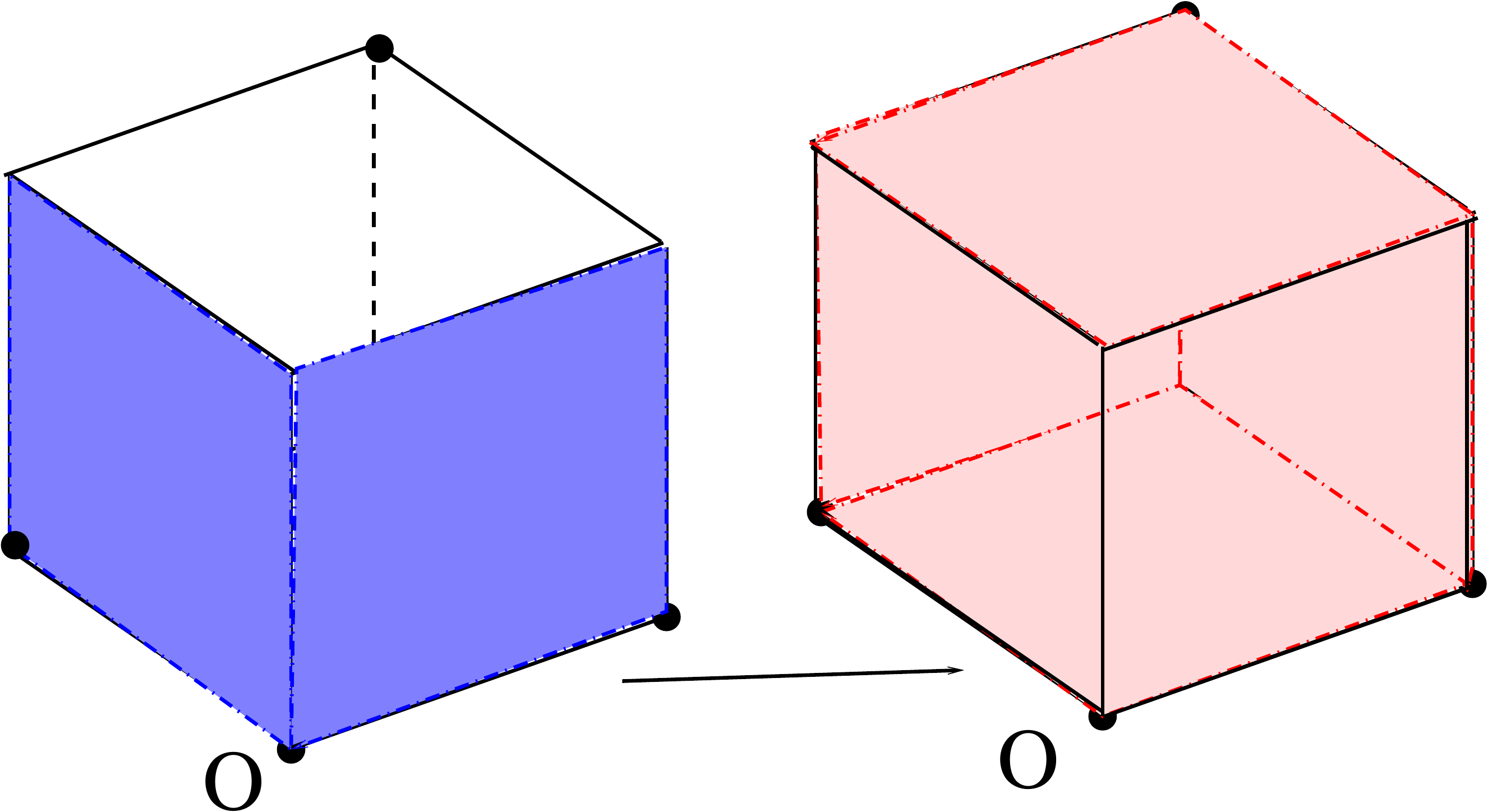}
\end{center}
\caption{\sl $N$ is ambient isotopic to $N_1$.} 
\label{F7}
\end{figure}

\noindent Clearly our gridded surface $N$ is ambient isotopic to a gridded surface $N_1=(N\setminus (F_{a,b}\cup F_{a,c}))\cup (F_{b,c}\cup F'_{a,b}\cup F'_{a,c}\cup F'_{b,c}))$ (see Figure \ref{F7}).

\noindent Observe that ${\cal{F}}(N_1)=({\cal{F}}(N)\setminus (F_{a,b}\cup F_{a,c}))\cup F_{b,c}$. $\square$\\

\noindent By the previous claim, our squared-star $\cal{F}(N)$ described by  $F_{a,b}\rightarrow F_{a,c}\rightarrow F_{b,c}$ is isotopic to 
$F_{a,b}\rightarrow F_{a,c}\rightarrow F_{-a,b}\rightarrow F_{-a,c}$; \emph{i.e.} 
\begin{equation}\label{1}
\begin{split}
F_{a,b}\rightarrow F_{a,c}\rightarrow F_{b,c}&\overset{c}\sim F_{a,b}\rightarrow F_{a,c}\rightarrow F_{-a,b}\rightarrow F_{-a,c}\\
&\overset{c}\sim F_{a,b}\rightarrow F_{a,-b}\rightarrow F_{-b,c}\rightarrow F_{-a,b}\rightarrow F_{-a,c}\\
&\overset{c}\sim F_{a,b}\rightarrow F_{a,-b}\rightarrow F_{-a,b}\rightarrow 
F_{-a,-b}.\\
\end{split}
\end{equation}

\noindent If $n=4$, then  ${\cal{F}}$ is described as a path as follows
\begin{equation}\label{2}
 F_{a,b}\rightarrow F_{a,c}\rightarrow F_{b,d} \rightarrow F_{c,d}.
\end{equation}
\noindent Then applying the above claim, we have that \eqref{2} is isotopic to
$$
F_{a,b}\rightarrow F_{a,c}\rightarrow F_{b,c}, 
$$
which by \eqref{1}, is isotopic to 
$F_{a,b}\,\rightarrow \, F_{a,-b}\,\rightarrow \, F_{-a,b}\,\rightarrow \, F_{-a,-b}$.\\

\noindent We assume that the Theorem holds for any squared-star consisting on $k$ faces, $k\leq n-1$. Suppose that $\cal{F}(N)$ consists of $n$ faces, \emph{i.e.}
\begin{equation}\label{3}
F_{a_1,b_1}\rightarrow F_{a_2,b_2}\rightarrow\, \ldots \, \rightarrow F_{a_{n-1},b_{n-1}} \rightarrow F_{a_n,b_n}.
\end{equation}
Since $N$ is a manifold, there exists $a_i$, $1\leq i\leq n-1$ such that $a_i=a_n$. For simplicity we will assume that $i=n-1$,
hence the squared-star \eqref{3} is
equal to
\begin{equation}\label{4}
F_{a_1,b_1}\rightarrow F_{a_2,b_2}\rightarrow\, \ldots \, \rightarrow F_{a_{n},b_{n-1}} \rightarrow F_{a_n,b_n}.
\end{equation}
\noindent By the above claim, the squared-star \eqref{4} is isotopic to
\begin{equation}\label{5}
F_{a_1,b_1}\rightarrow F_{a_2,b_2}\rightarrow\, \ldots \, \rightarrow F_{b_{n},b_{n-1}},
\end{equation}
\noindent and applying the induction hypothesis on \eqref{4} and \eqref{1}, we have that $\cal{F}(N)$ is isotopic to 
$F_{a,b}\,\rightarrow \, F_{a,-b}\,\rightarrow \, F_{-a,b}\,\rightarrow \, F_{-a,-b}$. Therefore, $slk(N)$ is isotopic to the boundary of the square $[-1,1]^2$ on 
 $ab$-plane, hence $slk(N)$ is unknotted simple closed curve. $\square$

\begin{coro}
${\cal{F}}(N)$ is topologically locally flat. $\square$
\end{coro}

\begin{theo}\label{flat}
Let $N$ be a gridded surface such that $0\in N$. Then $N$ is Whitehead locally flat at $0$.
\end{theo}
\noindent{\it Proof.} Consider the squared-link $slk(N)$ of $0$. By the Theorem \ref{unknotted}, we know that $slk(N)$ is an unknotted simple closed curve, then 
there exists  a smooth closed curve $J$ isotopic to $slk(N)$ such that $J$ is 
${\cal{C}}^0$-arbitrarily close to $slk(N)$. This is because
we can round the corners at the vertices of $slk(N)$ in an arbitrarily small neighborhoods of them (see \cite{douady}). Hence
$J$ is the boundary of a smooth disk $U_0$, ${\cal{C}}^0$-arbitrarily close to $\Int ({\cal{F}}(N))$. Let $\psi: U_0\times [0,1]\rightarrow  \Int ({\cal{F}}(N))$ be this isotopy. Then by \cite{whitehead} there exists a locally transverse 2-field 
 ${\text{{\bf v}}}(x)$ to each point of $x\in U_0$, hence  ${\text{{\bf v}}}(x)$ is also transverse to $y=\psi (x,1)\in \Int ({\cal{F}}(N))$. Therfore, the
 result follows. $\square$\\

\noindent We are ready to prove Theorem 1.

\begin{main}\label{main}
Any closed, oriented, gridded surface $N$ in $\mathbb{R}^{4}$ is smoothable. More precisely $N$
 admits a transverse field of 2-planes and therefore by a theorem of J. H. C. Whitehead there is an arbitrarily small topological isotopy that moves $N$ 
 onto a smooth surface in $\R^4$ (see  \cite{pugh}, \cite{whitehead}).
\end{main}
\noindent{\it Proof.} We will prove that $N$ admits a  global transverse 2-field  and   by  \cite{whitehead} it will follow that $N$ is smoothable. 
We will construct a local transverse 2-field at each point $x\in N$ in such a way that we can define a global transverse 2-field at $N$.  For the sake of simplicity, we divide each square $F\subset N$ of $\cal{C}$ into $m^2$ squares $S_i$, $i=1,2,\ldots ,m^2$; \emph{i.e.} each 
edge of $N$ is subdivided into $m$ equal segments. Let $x\in N$.\\

\noindent Case 1. The point $x$ lies on some square $S\subset\Int (F)$.  Consider the plane $P$ parallel to the support plane of $F$.
Then the map ${\text{{\bf v}}}_S:S\rightarrow G(2,4)$ defined as $y\mapsto y+P^{\bot}$, $y\in S$, where $P^{\bot}$ is the orthogonal plane to $P$. Thus
${\text{{\bf v}}}_S(y)$ is a local transverse 2-field at $y\in \Int(S)$ and 
for $y\in\partial S$, we have that it is a plane orthogonal to the corresponding square face.\\

\noindent Case 2. The point $x$ lies on  some square $S_1\subset F_1$ such that $S_1$ intersects some edge $e\subset N\subset {\cal{C}}$ and
none point of $S_1$ is a vertex of $\cal{C}$. We have that $e$ is the intersection of
 two faces $F_1$ and $F_2$ of $N$. Let $S_2\subset F_2$ be a square of $N$ such that $S_1\cap S_2\subset e$. Consider $l_1$ and $l_2$ edges of $S_1$ and $S_2$, respectively; such that
 $l_1\cap l_2=z\in e$; so $l_1\times l_2$ is a square. Let $C$ be the circle of radius $\frac{1}{m}$ and centered at the opposited vertex $z'$ to $z$ in $l_1\times l_2$, then
$C_{1,2}\subset C\cap l_1\times l_2$ is an arc such that, by construction, $l_1$ and $l_2$ are tangent to it. Let $r_x=x-z'$ be a vector, where $x\in l_1\cup l_2$, and 
let $e'$, $l_1'$ and $l_2'$ be vectors parallel to $e$, $l_1$ and $l_2$ respectively, at the origin. Observe that these three vectors are canonical vectors up to scale. Take the remaining canonical vector
$v_{1,2}$ and consider the plane $P_x=\langle r_x,v_{1,2}\rangle$. In general, for any  $y\in S_1\cup S_2$, there exists $x\in l_1\cup l_2$ such that $y-x$ is parallel to $e$, hence $P_y=\langle r_x,v_{1,2}\rangle$.
We define ${\text{{\bf v}}}_{S_1\cup S_2}:S_1\cup S_2\rightarrow G(2,4)$ given by $y\mapsto y+P_y$, $y\in S_1\cup S_2$. Then by construction, it is a local transverse 2-field at $y\in\Int (S_1\cup S_2)$ and
for $y\in\partial (S_1\cup S_2)\setminus \{l_1,l_1+\frac{1}{m}e,l_2, l_2+\frac{1}{m}e\}$, we have that ${\text{{\bf v}}}_{S_1\cup S_2}(y)$ is a plane orthogonal to the corresponding square face.\\

\begin{figure}[h]  
\begin{center}
\includegraphics[height=3cm]{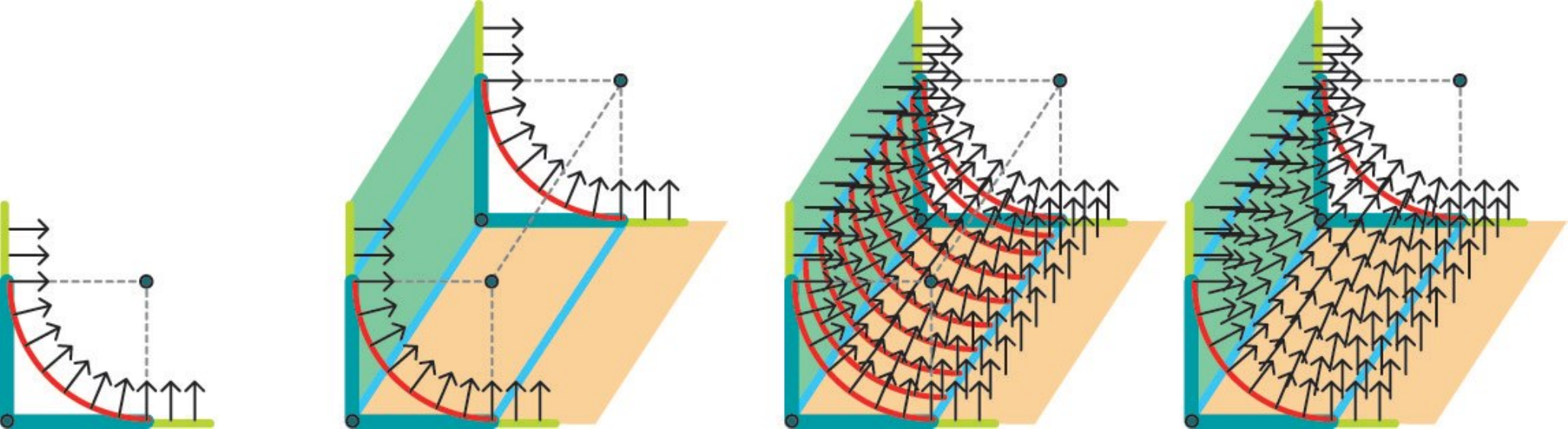}
\end{center}
\caption{\sl Vector field on the edge.} 
\label{F2}
\end{figure}

\noindent Case 3. The point $x$ is a vertex on the canonical cubulation $\cal C$, then $x=\cap_{j=1}^{r} F_{i_j}$, where $F_{i_j}$ is a 2-face of $N$.
We will suppose for simplicity that $x$ is the origin $0=(0,0,0,0)\in\mathbb{R}^4$.  As before, 
we will consider that  $x=\cap_{j=1}^{r} S_{j}$, where $S_{j}\subset F_{i_j}$ is a square; so ${\cal{F}}_m(N)=\cup_{j=1}^{r} S_{j}$. Let $e_j=S_{j}\cap S_{j+1}$ ($r+1=1$) be an edge whose end point is $x_j$. Then $slk_m(N)=\partial {\cal{F}}_m(N)$ consists on the union of a finite number of  edges; in fact, two edges for each square $S_j$. More precisely, $slk_m(N)=(\cup_{j=1}^{r} l_{j_1}\cup (\cup_{j=1}^{r} l_{j_2})$,
where $l_{j_1},\,l_{j_2}\subset S_j $,  $l_{j_1}$ is parallel to $e_{j+1}$ and $l_{j_2}$ is parallel to $e_{j}$ ($r+1=1$); in particular $x_j=l_{j_2}\cap l_{{j+1}_1}$.\\
Notice that using the same argument of case 2, we can define  ${\text{{\bf v}}}_{l_{j_2}\cup l_{{j+1}_1}}:l_{j_2}\cup l_{{j+1}_1}\rightarrow G(2,4)$ such that 
it is a local transverse 2-field at $y\in\Int (l_{j_2}\cup l_{{j+1}_1})$ and
for $y\in\partial (l_{j_2}\cup l_{{j+1}_1})$, we have that ${\text{{\bf v}}}_{l_{j_2}\cap l_{{j+1}_1}}(y)$ is a plane orthogonal to the corresponding square face.\\

\noindent By the above, we can define ${\text{{\bf v}}}_{slk_m(N)}:slk_m(N)\rightarrow G(2,4)$ such that if $x\in l_{j_2}\cup l_{{j+1}_1}$ then
${\text{{\bf v}}}_{slk_m(N)}(x):={\text{{\bf v}}}_{l_{j_2}\cup l_{{j+1}_1}}(x)$. 
Observe that by construction, ${\text{{\bf v}}}_{slk_m(N)}$ is well-defined continuous transverse 2-field.\\
By Theorem \ref{flat}, we know that there exists a local transverse 2-vector field ${\text{{\bf v}}}_{{\cal{F}}_m(N)}:{\cal{F}}(N)\rightarrow G(2,4)$, such that 
it is a local transverse 2-field at $y\in\Int ({\cal{F}}_m(N))$ and
for $y\in\partial (S_1\cup S_2)$, we have that ${\text{{\bf v}}}_{{\cal{F}}_m(N)}(y)= {\text{{\bf v}}}_{slk_m(N)}(y)$.\\

\noindent Next, we define $\frak{V}: N\rightarrow G(2,4)$ as $\frak{V}(y):= {\text{{\bf v}}}_{U}(y)$ if $y\in U$, where $U$ can be either a square $S$ (case 1), two squares $S_1\cup S_2$ (case 2) or a squared-star ${\cal{F}}_m(N)$ (case 3). Notice that by construction $\frak{V}(x)$
is a well-defined global transverse 2-field at $N$. The continuity of  $\frak{V}(x)$ follows from the fact that  by construction $\frak{V}(x)\rightarrow \frak{V}(y)$ as
 $x\rightarrow y$. Therefore $N$ is smoothable. $\square$

\noindent J. P. D\'iaz. {\tt Instituto de Matem\'aticas, Unidad Cuernavaca}. Universidad Nacional Au\-t\'o\-no\-ma de M\'exico.
Av. Universidad s/n, Col. Lomas de Chamilpa. Cuernavaca, Morelos, M\'exico, 62209.

\noindent {\it E-mail address:} juanpablo@matcuer.unam.mx\\

\noindent G. Hinojosa. {\tt Centro de Investigaci\'on en Ciencias, IICBA}. Universidad Aut\'onoma del Estado de Morelos. Av. Universidad 1001, Col. Chamilpa.
Cuernavaca, Morelos, M\'exico, 62209. 

\noindent {\it E-mail address:} gabriela@uaem.mx\\

\noindent R. Valdez. {\tt Centro de Investigaci\'on en Ciencias, IICBA}. Universidad Aut\'onoma del Estado de Morelos. Av. Universidad 1001, Col. Chamilpa.
Cuernavaca, Morelos, M\'exico, 62209. 

\noindent {\it E-mail address:} valdez@uaem.mx\\

\noindent A. Verjovsky. {\tt Instituto de Matem\'aticas, Unidad Cuernavaca}. Universidad Nacional Au\-t\'o\-no\-ma de M\'exico.
Av. Universidad s/n, Col. Lomas de Chamilpa. Cuernavaca, Morelos, M\'exico, 62209.

\noindent {\it E-mail address:} alberto@matcuer.unam.mx
\end{document}